\documentclass[10pt]{amsart}
\usepackage{amsmath}
\usepackage{amssymb}
\usepackage{amsfonts}
\usepackage{amsthm}
\usepackage{enumerate}
\usepackage[all]{xy}
\usepackage[latin1]{inputenc}
\usepackage{graphicx}
\usepackage{latexsym}
\usepackage{mathdots}
\input xy

\newcommand{\rep}{{\rm rep}}

\newcommand{\Hom}{{\rm Hom}}
\newcommand{\Ext}{{\rm Ext}}
\newcommand{\End}{{\rm End}}

\newcommand{\N}{\mathbb{N}}
\newcommand{\Z}{\mathbb{Z}}
\newcommand{\C}{\mathcal{C}}
\newcommand{\A}{\mathbb{A}}

\def\id{\hbox{1\hskip -3pt {\sf I}}}

\def\Ga{\hbox{$\mathit\Gamma$}}
\def\Da{\hbox{${\mathit\Delta}$}}
\def\Sa{\hbox{${\mathit\Sigma}$}}
\def\Oa{\hbox{${\mathit\Omega}$}}

\def\GaC{{\mathit\Gamma}_{_{\hspace{-1pt}\mathcal{C}}}}

\makeatletter
    \newtheoremstyle{definition}
        {5pt}
        {3pt}
        {}
        {0pt}
        {\scshape}
        {.}
        {5pt}
        {\thmname{#1} \thmnumber{#2} \thmnote{[#3]}} 

\newtheoremstyle{theorems}
        {5pt}
        {3pt}
        {\itshape}
        {0pt}
        {\scshape}
        {.}
        {5pt}
        { \thmname{#1} \thmnumber{#2}\thmnote{[#3]}} 

\swapnumbers \theoremstyle{theorems}
\newtheorem{Theo}{Theorem}[section]

\newtheorem{Cor}[Theo]{Corollary}
\newtheorem{Lemma}[Theo]{Lemma}

\theoremstyle{definition}

\newtheorem{Defn}[Theo]{Definition}

\begin{document}

\title[Standard components]{\sc Standard components of a Krull-Schmidt category}

\vspace{-5pt}

\author[Shiping Liu]{Shiping Liu}
\author[Charles Paquette]{Charles Paquette}


\thanks{The first named author is supported in part by the NSERC of Canada, while the second named author is supported in full by NSERC and AARMS}


\maketitle

\vspace{-23pt}

\begin{abstract}

We provide criteria for an Auslander-Reiten component having sections of a Krull-Schmidt category to be standard. Specializing to the cate\-gory of finitely presented representations of a strongly locally finite quiver and its bounded derived category, we obtain many new types of standard Auslander-Reiten components. An application to the module category of a finite-dimensional algebra yields some interesting results.

\end{abstract}

\medskip

\section*{Introduction}

Standard Auslander-Reiten components of the module category  of a finite dimensional algebra are extremely interesting, since the maps between modules in such a component can be described in a simple combinatorial way; see \cite{BoG, R3}. This kind of components appears mainly for representation-finite algebras, hereditary algebras, tubular algebras and tilted algebras; see \cite{R3}, and each of them has at most finitely many non-periodic Auslander-Reiten orbits; see \cite{Skow}. In particular, the regular ones are stable tubes or of shape $\Z\Da$ with $\Da$ a finite acyclic quiver.

\medskip

On the other hand, the Auslander-Reiten theory has been extended to Krull-Schmidt categories; see \cite{Bau, L3}. It is natural to expect that new types of standard Auslander-Reiten components will appear in this context.
Indeed, in the most gene\-ral setup, we shall find various criteria for such an Auslander-Reiten component having sections to be standard. In particular, an Auslander-Reiten component which is a wing or of shape $\N\A^+_\infty$, $\N^-\A^-_\infty$ or $\Z\A_\infty$ is standard if and only if its quasi-simple objects are pairwise orthogonal bricks. Specializing to $\rep^+(Q)$, the category of finitely presented representations of a connected strongly locally finite quiver $Q$, we prove that the preprojective component and the preinjective components are standard; and every component is standard in case $Q$ is of finite or infinite Dynkin type. Applying this to the bounded derived category $D^b(\rep^+(Q))$ of $\rep^+(Q)$, we show that the connecting component is standard; and every component is standard in the Dynkin case. All these particularly establish the existence of standard  Auslander-Reiten components which are wings or of shapes $\N\A^+_\infty$, $\N^-\A^-_\infty$ and $\Z\Da$ with $\Da$ any strongly locally finite quiver without infinite paths. Furthermore, specialized to the module category ${\rm mod}\hspace{0.4pt}A$ of a finite dimensional algebra $A$, our criteria become surprisingly nice and easy to verify; see (\ref{maincor}). As a consequence, an Auslander-Reiten component with sections of ${\rm mod}\hspace{0.4pt}A$ is standard if and only if it is generalized standard, if and only if it is the connecting component of a tilted factor algebra of $A$. Finally, we remark that some of our results will be applied in the future to study cluster categories of infinite Dynkin types.

\section{Standard components having sections}

\smallskip

Throughout this paper, $k$ stands for an arbitrary field. A $k$-category is a category in which the morphism sets are $k$-vector spaces and the composition of morphisms is $k$-bilinear. A $k$-category is called {\it Hom-finite} if its morphism spaces are all finite dimensional over $k$, and {\it Krull-Schmidt} if every non-zero object is a finite direct sum of objects with a local endomorphism algebra.

\medskip

For the rest of this section, let $\C$ stand for a Hom-finite Krull-Schmidt additive $k$-category. The {\it radical} morphisms in $\C$ are those in the Jacobson radical ${\rm rad}(\C).$ One calls ${\rm rad}^\infty(\C)=\cap_{n\ge 1}\,{\rm rad}^n(\C)\vspace{1pt}$ the {\it infinite radical} of $\C$, where ${\rm rad}^n(\C)$ is the $n$-th power of ${\rm rad}(\C)$. Two objects $X, Y$ in $\C$ are said to be \emph{orthogonal} if $\Hom_{\hspace{0.4pt}\mathcal{C}}(X, Y) = 0$ and $\Hom_{\hspace{0.4pt}\mathcal{C}}(Y, X)=0$. If $X\in \C$ is indecomposable, then the division algebra $k_{\hspace{-0.5pt}_X}\hspace{-2pt} \vspace{1pt} ={\rm End}(X)/ \hspace{0.5pt}{\rm rad}(X, X)$ is called the {\it automorphism field} of $X$, and we shall call $X$ a {\it brick} provided that $\End_{\mathcal{C}}(X)$ is trivial, that is, $\End_{\mathcal{C}}(X)\cong k.$
Let $f: X \to Y$ be a morphism in $\C$. One says that $f$ is {\it irreducible} if it is neither a section nor a retraction, and any factorization $f=gh$ implies that $h$ is a section or $g$ is a retraction.
Moreover, $f$ is called \emph{left almost split} if it is not a section and every non-section morphism $g: X \to M$ in $\C$ factors through $f$; \emph{left minimal} if every endomorphism $h$ of $Y$ such that
$f=hf$ is an automorphism. In a dual manner, one defines $f$ to be {\it right almost split} and {\it right minimal.} Further,
$f$ is called a \emph{source morphism for $X$} if it is left minimal and left almost split, and a {\it sink morphism for $Y$} if it is right minimal and right almost split. A sequence of morphisms $$\xymatrix{X \ar[r]^f & Y \ar[r]^g & Z}$$ in $\C$ with $Y\ne 0$ is called {\it almost split} provided that $f$ is a source morphism and a pseudo-kernel of $g$, while $g$ is a sink morphism and a pseudo-cokernel of $f$; see \cite[(1.3)]{L3}.
In case $\C$ is abelian or triangulated, the definition of an almost split sequence given here coincides somehow with the classical one; see \cite[(1.5), (6.1)]{L3}.

\medskip

We shall make a frequent use of the following easy result.

\smallskip

\begin{Lemma}\label{auto-field} Let  $\C$ have  an almost split sequence as follows$\,:$
$$\xymatrixcolsep{25pt}\xymatrix{X \ar[r]^-{f_1 \choose f_2} & Y_1\amalg Y_2 \ar[r]^-{(g_1, g_2)} & Z.}\vspace{-2pt}$$

\begin{enumerate}[$(1)$]

\item There exists a $k$-linear isomorphism $k_{\hspace{-0.8pt}_X}\hspace{-1pt}\cong k_{\hspace{-1pt}_Z}\vspace{0.5pt}$.

\item If $u: M\to Y_1$ is a morphism in $\C$ such that $g_1 u=0$, then there exists some $w: M\to X$ such that $u=f_1w$ and $f_2w=0$.

\item If $v: Y_1\to N$ is a morphism in $\C$  such that $v f_1=0$, then there exists some $w: Z\to N$ such that $v=wg_1$ and $wg_2=0.$

\end{enumerate}
\end{Lemma}

\noindent{\it Proof.} Statement (1) is implicitly stated and proved in the proof of \cite[(2.1)]{L3}.\vspace{1pt} Let $u: M\to Y_1$ be such that $g_1 u=0$. Then $(g_1, g_2) {u\choose 0}=0,$ and hence there exists some $w: M\to X$ such that ${u\choose 0}={f_1\choose f_2}w.\vspace{1pt}$ This proves Statement (2). Dually, we can show Statement (3). The proof of the lemma is completed.

\medskip

The \emph{Auslander-Reiten quiver} $\Ga_{_{\hspace{-1pt}\mathcal{C}}}$ of $\C$ is first defined to be a valued translation quiver as follows. The vertex set is a complete set of the representatives of the isomorphism
classes of the indecomposable objects in $\mathcal C$. For
vertices $X$ and $Y$, we write ${d\hspace{0.4pt}'}\hspace{-2pt}_{_{XY}}\vspace{0.5pt}$ and $d_{_{XY}}\vspace{0.5pt}$ for the dimensions of $${\rm irr}(X, Y)={\rm rad}(X, Y)/\hspace{0.5pt}{\rm rad}^2(X, Y)\vspace{0.5pt}$$ over $k_{\hspace{-0.5pt}_X}\hspace{-1pt}$ and $k_{_Y}\hspace{-1pt}$ respectively, and draw a unique valued arrow $X\to Y$ with
valuation $(d_{_{XY}}, {d\hspace{0.4pt}'}\hspace{-2pt}_{_{XY}}\vspace{0.5pt})$ if and only if
$d_{_{XY}}>0.$ The translation $\tau$ is defined so that
$\tau Z=X$ if and only if $\mathcal C$ has an almost split
sequence $\xymatrixcolsep{15pt}\xymatrix{X\ar[r] & Y \ar[r] & Z.}$ A valuation $(d_{_{XY}},
{d\hspace{0.4pt}'}\hspace{-2pt}_{_{XY}}\vspace{0.5pt})$ is called {\it symmetric} if
$d_{_{XY}}={d\hspace{0.4pt}'}\hspace{-2pt}_{_{XY}}\vspace{0.5pt}$, and {\it trivial} if
$d_{_{XY}}={d\hspace{0.4pt}'}\hspace{-2pt}_{_{XY}}\vspace{0.5pt}=1$. Next, $\GaC$ is modified in such a way that each symmetrically valued arrow $X\to Y$ is replaced by $d_{_{XY}}$ unvalued arrows from $X$ to $Y$. That is,
$\Ga_{_{\hspace{-1pt}\mathcal{C}}}$ becomes a partially valued translation quiver in which all valuations are non-symmetric; see \cite[(2.1)]{L3}.

\medskip

Let $\Sa$ be a convex subquiver of $\Ga_{_{\hspace{-1pt}\mathcal{C}}}$ in which every object has a trivial automorphism field. In particular, $d_{_{XY}}={d\hspace{0.4pt}'}\hspace{-2pt}_{_{XY}}\vspace{0.5pt}$ for all $X, Y\in \Sa$. By our construction, $\Sa$ is a non-valued translation quiver with possible multiple arrows. Thus, one can define the {\it path category} $k\Sa$ and the {\it mesh category} $k(\Sa)$ of $\Sa$ over $k\hspace{0.6pt};$ see, for example, \cite[(2.1)]{R3}. In the sequel, for $u\in k\Sa$, we shall write $\overline{u}$ for its image in $k(\Sa)$.

\medskip

\begin{Defn}\label{standardness}

Let $\Sa$ be a convex subquiver of $\GaC$, and let $\C(\Sa)$ be the full subcategory of $\C$ generated by the objects in $\Sa$. We shall say that $\Sa$ is {\it standard} provided that every object in $\Sa$ has a trivial automorphism field and there exists a $k$-equivalence $F: k({\Sa}) \stackrel{\sim}{\rightarrow}\C(\Sa)$, which acts identically on the objects.

\end{Defn}

\medskip

\begin{Lemma}\label{irr-basis}

Let $\Sa$ be a convex subquiver of $\GaC$, and let $F: k({\Sa}) \stackrel{\sim}{\rightarrow}\C(\Sa)$ be a $k$-equivalence acting identically on the objects. If $X, Y\in \Sa$, then the classes $F(\overline{\alpha})+{\rm rad}^2\hspace{-0.4pt}(X, Y)$ form a $k$-basis of $\,{\rm irr}(X, Y)$, where $\alpha$ ranges over the set of arrows from $X$ to $Y$.

\end{Lemma}

\noindent{\it Proof.} Let $X, Y\in \Sa$. For $1\le i\le 2,$ consider the $k$-subspace ${\mathcal I}^{\,(i)}(X, Y)$ of $k(\Sa)(X, Y)$ generated by the $\overline{p}$, where $p$ ranges over the set of paths of length $\ge i$ from $X$ to $Y$. Write $\Sa_1(X, Y)$ for the set of arrows from $X$ to $Y$. Since the mesh relations are sums of paths of length two, the classes $\overline\alpha + \mathcal{I}^{\,(2)}(X, Y)$, with $\alpha\in \Sa_1(X, Y),$ are $k$-linearly independent, and hence, they form a $k$-basis for $\mathcal{I}^{\,(1)}(X, Y)/\mathcal{I}^{\,(2)}(X, Y)$. Thus, $\mathcal{I}^{\,(1)}(X, Y)/\mathcal{I}^{\,(2)}(X, Y)$ and ${\rm irr}(X, Y)$ are of the same $k$-dimension.
Since $F$ induces a $k$-isomorphism $F: k(\Sa)(X, Y)\to \Hom_{\,\mathcal{C}}(X, Y)$, it is easy to see that $F$ induces a $k$-epimorphism $F: \mathcal{I}^{\,(1)}(X, Y)\to {\rm rad}(X, Y).\vspace{1pt}$ In particular, $F$ maps $\mathcal{I}^{\,(2)}(X, Y)$ into ${\rm rad}^2(X, Y)$. This yields a $k$-epimorphism

\vspace{-14pt}
$$\overline{F}: \mathcal{I}^{\,(1)}(X, Y)/\mathcal{I}^{\,(2)}(X, Y) \to {\rm irr}(X, Y): u + \mathcal{I}^{\,(2)}(X, Y) \mapsto F(u) +{\rm rad}^2(X, Y),\vspace{-4pt}$$ which is necessarily an isomorphism. The proof of the lemma is completed.

\medskip

Given a quiver $\Sa$ with no oriented cycle, one constructs a stable translation quiver $\Z \Sa$; see, for example,
\cite[(2.1)]{R3}. We denote by $\N\Sa$ the full translation subquiver of $\Z\Sa$ generated by the vertices $(n, x)$ with $n\ge 0$ and $x\in \Sa$, and by $\N^-\hspace{-3pt}\Sa$ the one generated by the vertices $(n, x)$ with $n\le 0$ and $x\in \Sa$. Now, let $\Ga$ be a connected component of $\GaC$. A connected full subquiver $\Da$ of $\Ga$ is called a \emph{section} if it is convex in $\Ga$, contains no oriented cycle, and meets every $\tau$-orbit in $\Ga$
exactly once. In this case, every object in $\Ga$ is uniquely written as $\tau^nX$ with $n\in \Z$ and $X\in \Da$, and
there exists a translation-quiver embedding $\Ga\to \Z\Da: \tau^nX\mapsto (-n, X);$ see \cite[(2.3)]{L4}. We denote by $\Da^-$ the full subquiver of $\Ga$ generated by the vertices $\tau^nX$ with $n>0$ and $X\in \Da$, and by $\Da^+$ the one generated by the vertices $\tau^nX$ with $n<0$ and $X\in \Da$. One says that $\Da$ is a \emph{right-most section} if $\Da^+=\emptyset;$ and \emph{left-most section} if $\Da^-=\emptyset.$

\medskip

In order to state and prove the following main result of this section, we need some terminology and notation. Firstly, an infinite path in a quiver is called {\it left infinite} if it has no starting point; and {\it right infinite} if it has no ending point. Secondly, given two (possibly empty) subquivers $\Sa, \Oa$ of $\GaC$, we shall write $\Hom_{\hspace{0.5pt}\mathcal{C}}(\Sa, \Oa)=0$ in case $\Hom_{\hspace{0.5pt}\mathcal{C}}(X, Y)=0$ for all possible objects $X\in \Sa$ and $Y\in \Oa$.

\medskip

\begin{Theo} \label{mainprop}

Let $\C$ be a Hom-finite Krull-Schmidt additive $k$-category, and let $\Ga$ be a connected component of $\GaC$ having a section $\Da.$ If $\Da^+$ has no left infinite path and $\Da^-$ has no right infinite path, then ${\it \Gamma}$ is standard if and only if
$\it\Delta$ is standard such that $\Hom_{\hspace{0.5pt}\mathcal{C}}(\Da^+, \Da\cup \Da^-)=0$ and $\Hom_{\hspace{0.5pt}\mathcal{C}}(\Da, \Da^-) = 0.$

\end{Theo}

\noindent{\it Proof.} Suppose that $\Da^+$ has no left infinite path and $\Da^-$ has no right infinite path. Assume first that $\Ga$ is standard. In particular, $\Da$ is standard. Since $\Ga$ embeds in $\Z\Da$, we see that $\Ga$ has no path from $X$ to $Y$ in case $X\in \Da^+$ and $Y\in \Da\cup\Da^-$, or $X\in \Da$ and $Y\in \Da^-$. This shows the necessity.

Assume conversely that $\mathit\Delta$ is standard such that $\Hom_{\,\mathcal{C}}(\Da^+,  \Da\cup\Da^-) =0$ and $\Hom_{\,\mathcal{C}}(\Da, \Da^-) = 0.$ In particular, every object in $\Da$ has a trivial endomorphism algebra. Being of the form $\tau^nX$ with $n\in \Z$ and $X\in \Da$, by Lemma \ref{auto-field}(1), every object in $\Ga$ has a trivial automorphism field. Since every object in $\Da^+$ has a sink morphism and $\Da$ is a section, $\Da^+$ is locally finite. By K\"{o}nig's Lemma, $\Da^+$ has only finitely many paths ending in any pre-fixed object. Thus, for each object $M\in \Da\cup \Da^+$, we may define an integer $n_{\hspace{-0.5pt}_M}\ge 0$ in such a way that $n_{\hspace{-0.5pt}_M}=0$ if $M\in \Da$; and otherwise, $n_{\hspace{-0.5pt}_M}-1$ is the maximal length of the paths of $\Da^+$ which end in $M$. The following statement is evident.

\vspace{1pt}

(1) {\it Let $p: X\rightsquigarrow Y$ be a non-trivial path in $\Ga$. If $X\in \Da\cup \Da^+$, then $Y\in \Da\cup \Da^+$ with
$n_{\hspace{-0.5pt}_X}\le n_{\hspace{-0.5pt}_Y}$, and the equality occurs if and only if $X, Y\in \Da$.}

\vspace{1pt}

For each $n\ge 0$, denote by $\Ga^{\,n}$ the full subquiver of ${\it \Gamma}$ generated by the vertices $X \in \Da\cup \Da^+$ with $n_{\hspace{-0.5pt}_X} \le n\vspace{1pt}$, which is clearly convex in $\Ga$. Moreover, denote by $\Ga^{\,+}$ the union of the $\Ga^{\,n}$ with $n\ge 0$, that is, the full subquiver of $\Ga$ generated by $\Da^+\cup \Da$.
The following statement is an immediate consequence of Statement (1).

\vspace{1pt}

(2) {\it If $p: X\rightsquigarrow Y$ is a non-trivial path in $\Ga^{\,n+1}$ with $n\ge 0$, then $X\in \Ga^{\,n},$ and consequently, $p\not\in \Ga^{\,n}$ if and only if $Y\not\in \Ga^{\,n}.$}

\vspace{1.5pt}

Now, let $F^{\,0}: k(\Da) \stackrel{\sim}{\to} \C(\Delta)$ be a $k$-linear equivalence, acting identically on the objects. Since $\Da$ contains no mesh of $\Ga$, we have $k(\Da)=k\Da$. Assume that $n\ge 0$ and $F^{\,0}$ extends to a full $k$-linear functor $F^{\,n}: k\Ga^{\,n}\to \C(\Ga^{\,n})$, acting identically on the objects and having a kernel generated by the mesh relations. In order to extend $F^{\,n}$ to $k\Ga^{\,n+1}$, we shall need the following statement.

\vspace{1pt}

(3) {\it If $f: X\to Y$ is a non-zero radical morphism in $\C(\Ga^{\,n+1})$, then $\Ga^{\,n+1}$ has a non-trivial path from $X$ to $Y$}. To show this, we claim that there exists $M\in \Da$, which is a predecessor of $Y$, such that $\Hom_{\hspace{0.5pt}\mathcal{C}}(X, M)\ne 0$. Indeed, suppose that this claim was false. Then, $Y\in \Da^+$. Since $\Da$ is a section in $\Ga$, every immediate predecessor of an object $\Da^+$ lies in $\Da^+\cup \Da$. Since every object in $\Da^+$ admits a sink morphism in $\C$, by factorizing the radical morphism $f$, we get a left infinite path
$$\xymatrixcolsep{16pt}\xymatrix{\cdots \ar[r] & Y_i\ar[r] & Y_{i-1} \ar[r] & \cdots \ar[r] & Y_1\ar[r] & Y}\vspace{-1pt}$$ in $\Da^+$ such that $\Hom_{\hspace{0.5pt}\mathcal{C}}(X, Y_i)\ne 0$ for all $i>0$, contrary to the hypothesis on $\Da^+$. Thus, $\Da$ does contain an object $M$ as claimed. Since $\Hom_{\hspace{0.5pt}\mathcal{C}}(\Da^+, \Da)=0$, we have $X\in \Da$. Since $k\Da\cong \C(\Da)$, there exists a path in $\Da$ from $X$ to $M$. This yields a non-trivial path in $\Ga^{\,n+1}$ from $X$ to $Y$. Statement (3) is established.

Fix $Z\in {\it\Gamma}^{\,n+1}\backslash {\it\Gamma}^{\,n}$. Observe that $Z\in \Da^+$ and $\tau Z\in \Da^+\cup \Da$. Thus, $k\Ga^{\,n+1}$ has a mesh relation $\delta\hspace{-0.5pt}_Z=\sum_{i=1}^r \beta_i\alpha_i,$ where $\alpha_i: \tau Z\to Y_i$, $i=1, \ldots, r$, are the arrows starting in $\tau Z$ and $\beta_i: Y_i\to Z$, $i=1, \ldots, r$, are the arrows ending in $Z$.
%
%
By Statement (2), $\tau Z, Y_1, \ldots, Y_r\in \Ga^{\,n}.$  Since $\tau Z$ admits a source morphism in $\C$, it follows from Lemma \ref{irr-basis} that $f=(F^{\hspace{0.4pt}n}(\alpha_1), \ldots, F^{\hspace{0.4pt}n}(\alpha_r))^T: \tau Z\to Y_1\oplus \cdots \oplus Y_r$ is a source morphism, which embeds in an almost split sequence
\vspace{-4pt}$$(*) \quad \xymatrix{\tau Z \ar[r]^-{f}
&Y_1\oplus \cdots \oplus Y_r\hspace{-6pt}} \xymatrixcolsep{13pt}\xymatrix{\ar[rr]^-{(g_1, \ldots, g_r)}&& Z
}$$
in $\C$; see \cite[(1.4)]{L3}. Set $F^{\hspace{0.4pt}n+1}(Z)=Z$, $F^{\hspace{0.4pt}n+1}(\varepsilon_{_Z})=\id_{_Z},$ where $\varepsilon_{_Z}$ is the {\it trivial} path at $Z$, and $F^{\hspace{0.4pt}n+1}(\beta_i)=g_i,\vspace{1pt}$ for $i=1, \ldots, s$. In view of Statement (2), we have defined $F^{\hspace{0.4pt}n+1}$ on the vertices, the trivial paths, and the arrows in $\Ga^{\,n+1}$. In an evident manner, we may extend $F^{\hspace{0.4pt}n}$ to a $k$-functor $F^{\hspace{0.4pt}n+1}: k\Ga^{\,n+1}\to \C(\Ga^{\,n+1})$, acting identically on the objects.

Let $u: Y\to Z$ be a non-zero radical morphism in $\C(\Ga^{\,n+1})$. By Statement (3),
$\Ga^{\,n+1}$ has a non-trivial path from $Y$ to $Z$, and hence $Y\in \Ga^{\,n}$ by Statement (2). If $Z\in \Ga^{\,n}$ then, by the induction hypothesis, $u=F^{\,n}(\rho)$ for some morphism $\rho: Y\to Z$ in $k\Ga^{\,n}$. Otherwise, $Z$ is the ending term of an almost split sequence $(*)$ as stated above. Then $u=\sum_{i=1}^r\, g_iu_i$, with morphisms $u_i: Y\to Y_i$ in $\mathcal C$.
Since $Y_i\in \Ga^{\,n}$, there exists $\rho_i: Y\to Y_i$ in $k\Ga^{\,n}$ such that $u_i=F^{\,n}(\rho_i)$, for $i=1, \ldots, r$. This yields $u=F^{\,n+1}(\sum_{i=1}^r\, \beta_i\rho_i)\vspace{0.5pt}$. That is, $F^{\,n+1}$ is full.

Next we shall show, for $\theta\in k\Ga^{\,n+1},$ that $F^{\,n+1}(\theta)=0$ if and only if $\theta$ lies in
the mesh ideal of $k\Ga^{\,n+1}$. In view of the induction hypothesis, we may assume that $\theta$ is non-zero of the form $\theta: Y\to Z$ with $Z\in\Ga^{\,n+1} \backslash \Ga^{\,n}$. In particular, $\Ga^{\,n+1}$ has a non-trivial path from $Y$ to $Z$. By Statement (2), $Y\in \Ga^{\,n}$. Suppose first that $\theta$ lies in the mesh ideal of $k\Ga^{\,n+1}$. For simplicity, we may assume that $\theta=\zeta\, \delta \, \sigma,$ where $\sigma, \delta, \zeta\in k\Ga^{\,n+1}$ with $\delta$ a mesh relation. If $\zeta$ has as a non-zero summand a multiple of a non-trivial path, then $\delta\in k\Ga^{\,n}$ by Statement (2). Hence, $F^{\,n+1}(\theta)=0$ by the induction hypothesis. Otherwise, $\delta$ is the mesh relation $\delta_{\hspace{-0.5pt}Z}$ as stated above, and $\theta=(\sum_{i=1}^r\beta_i\alpha_i)\eta$, where $\eta: Y\to \tau Z$ is a morphism in $k\Ga^{\,n}$. Since $(*)$ is an almost split sequence, we obtain $F^{\,n+1}(\theta)=0$.

Suppose conversely that $F^{\,n+1}(\theta)=0$. Consider the mesh relation $\delta_{\hspace{-0.5pt}Z}$ and the almost split sequence $(*)$ as stated above. Then $\theta=\sum_{i=1}^r \beta_i\theta_i$, where $\theta_i: Y\to Y_i$ is in $k\Ga^{\,n}$. Since
$\sum_{i=1}^r \,F^{\,n+1}(\beta_i) F^{\,n}(\theta_i)=F^{\,n+1}(\theta)=0$, there exists $v: Y\to \tau Z$ in $\C$ such that $F^{\,n}(\theta_i)=F^{\,n}(\alpha_i)\,v$, for $i=1, \ldots, r$. Since $F^{\,n}$ is full, $v=F^{\,n}(\eta)$ with $\eta: Y\to \tau Z$ in $k\Ga^{\,n}$. Hence $F^{\,n}(\theta_i)=F^{\,n}(\alpha_i\eta)$, and by the induction hypothesis, $\theta_i-\alpha_i\eta$ lies in the mesh ideal of $k\Ga^{\,n}$, $i=1, \ldots, r$. As a consequence, 
$$\theta=\textstyle\sum_{i=1}^r\beta_i(\theta_i-\alpha_i\eta)+ (\textstyle\sum_{i=1}^r \beta_i\alpha_i)\eta$$ lies in the mesh ideal of $k\Ga^{\,n+1}.$ This shows that $F^{\,n+1}$ is full and its kernel is generated by the mesh relations. By induction, $F^{\,0}$ extends to a full $k$-functor $F^+: k\Ga^+\to \C(\Ga^+)$, acting identically on the objects and having a kernel generated by the mesh relations.

Finally, for each object $N\in\Ga$, we may define $m_{_N}\ge 0$ so that $m_{_N}=0$ if $N\in \Ga^{\,+}$; and otherwise, $m_{_N}-1$ is the maximal length of the paths of $\Da^-$ which start in $N$. For $m\ge 0$, denote by $\Ga^{\hspace{0.4pt}(\hspace{-0.5pt}m)}$ the full subquiver of $\Ga$ generated by the objects $Y$ with $m_{_Y}\le m.$
Then $\Ga$ is the union of the $\Ga^{\hspace{0.4pt}(\hspace{-0.5pt}m)}$ with $m\ge 0\,.$ In a dual manner, we may apply the induction on $m$ to show that $F^+$ extends to a full $k$-functor $F: k\Ga\to \C(\Ga),$ which acts identically on the objects and has a kernel generated by the mesh relations. The proof of the theorem is completed.

\medskip

%

\begin{Lemma}\label{no-inf-path} Let $\Ga$ be a connected component of $\GaC$, containing a section $\Da.$

\vspace{-2pt}

\begin{enumerate}[$(1)$]

\item If $\Da$ has no left infinite path, then $\Da^+$ has no left infinite path.

\item If $\Da$ has no right infinite path, then $\Da^-$ has no right infinite path.

\end{enumerate}

\end{Lemma}

\noindent{\it Proof.} It suffices to prove Statement (1). Suppose that $\Da^+$ has a left infinite path

\vspace{-8pt}

$$\xymatrixcolsep{14pt}\xymatrix{\cdots \ar[r] & \tau^{-n_i}X_i \ar[r] & \cdots \ar[r] & \tau^{-n_1}X_1\ar[r] & \tau^{-n_0}X_0,} \vspace{-1pt}$$ where
$X_i\in \Da$ and $n_i>0$. Since $\Ga$ embeds in $\Z \Da$; see \cite[(2.3)]{L4}, we see that $n_i\le n_{i-1}$ for all $i> 0$. As a consequence, there exists $r\ge 0$ such that $n_i=n_r$ for $i\ge r$. Thus $\xymatrixcolsep{16pt}\xymatrix{\cdots \ar[r] & X_i \ar[r] & \cdots \ar[r] & X_r}$ is a left infinite path in $\Da$.
The proof of the lemma is completed.

\medskip

We shall say that a sink morphism in $\C$ is {\it proper} if it either is a monomorphism or fits in an almost slit sequence; and dually, a source morphism is {\it proper} if it either is an epimorphism or fits in an almost slit sequence. Observe that sink or source morphisms in an abelian category are all propre. The following result is a generalization of Lemma 3 stated in \cite[(2.3)]{R3}.

\medskip

\begin{Theo}\label{Theorem2}

Let $\C$ be a Hom-finite Krull-Schmidt additive $k$-category. Let $\Ga$ be a connected component of $\GaC$, and let $\Da$ be a section of $\Ga$ in which every object has a trivial automorphism field and admits a proper sink morphism as well as a proper source morphism. If $\Da$ has no infinite path, then $\Ga$ is standard if and only if $\,\Hom_{\hspace{0.5pt}\mathcal{C}}(\Da^+, \Da^-)=0.$

\end{Theo}

\noindent{\it Proof.} Suppose that $\Da$ has no infinite path. By Lemma \ref{no-inf-path}, $\Da^+$ has no left infinite path and $\Da^-$ has no right infinite path. We shall need the following statement.

\vspace{1pt}

Sub-Lemma: {\it Let $M\in \Ga$ with $\Hom_{\hspace{0.5pt}\mathcal{C}}(M, \Da^-)=0,$ and let $N\in \Da$. If $\C$ has a non-zero radical morphism $f: M\to N$, then $\Ga$ has a non-trivial path $M \rightsquigarrow N$.}

\vspace{1pt}


Indeed, suppose that $\Ga$ has no non-trivial path from $M$ to $N$. By assumption, $N$ admits a sink morphism $g=(g_1, \cdots, g_r): N_1\oplus \cdots \oplus N_r\to N$, where $N_i\in \Ga$. If $f: M\to N$ is non-zero and radical, then $f=\sum_{\,i=1}^{\,r}\,g_if_i$, with $f_i: M\to N_i$ in $\C$. We may assume that $f_1$ is non-zero. Since $\Da$ is a section, $N_1\in \Da\cup \tau \hspace{-2pt}\Da$; see \cite[(2.2)]{L4}. Since $\Hom_{\hspace{0.5pt}\mathcal{C}}(M, \Da^-)=0,$ we have $N_1\in \Da$. Since $\Ga$ has no path from $M$ to $N_1$, we see that $f_1$ is radical. Repeating this process, we see that $\Da$ contains an infinite path ending in $N$, a contradiction. This proves the sub-lemma.

Now, assume that $\Hom_{\hspace{1pt}\mathcal{C}}(\Da^+, \Da^-)=0.$ We deduce from the above sub-lemma that $\Hom_{\,\mathcal{C}}(\Da^+, \Da )=0.$ Using the dual statement, we obtain $\Hom_{\hspace{0.5pt}\mathcal{C}}(\Da, \Da^-)=0.$
It remains to construct a $k$-linear equivalence $F: k\Da\to \C(\Da)$. Since the objects in $\Da$ have a trivial automorphism field, so do the objects in $\Ga$. Set $F(X)=X$ and $F(\varepsilon_{\hspace{-0.5pt}_X})=\id_{\hspace{-0.5pt}_X}\vspace{1pt}$ for $X\in \Da$. Let $X, Y\in \Da$ with $d=d_{\hspace{-1pt}_{XY}}>0$. If $\alpha_i: X\to Y$, $i=1, \ldots, d$, are the arrows from $X$ to $Y$, then we choose irreducible morphisms $f_{\alpha_i}: X\to Y$ such that $f_{\alpha_1}+{\rm rad}^2(X, Y), \ldots, f_{\alpha_r}+{\rm rad}^2(X, Y)$ form a $k$-basis of ${\rm irr}(X, Y)$, and set $F(\alpha_i)=f_{\alpha_i}$, $i=1, \ldots, d$. In an evident manner, we obtain a $k$-linear functor $F: k\Da\to \C(\Da)$.

We claim that $F$ induces a $k$-isomorphism $F_{_{XY}}: \Hom_{\hspace{0.4pt}k\it\Delta}(X, Y)\to \Hom_{\,\mathcal{C}}(X, Y),$ for any $X, Y\in \Da$. Since every object in $\Da$ admits a sink morphism and a source morphism, $\Da$ is locally finite. Having no infinite path, by K\"onig's Lemma, $\Da$ has at most finitely many paths from $X$ to $Y$. Define an integer $n_{\hspace{-0.5pt}_{XY}}$ in such a way that $n_{\hspace{-0.5pt}_{XY}}=-1$ if $\Da$ has no path from $X$ to $Y$; and otherwise, $n_{\hspace{-0.5pt}_{XY}}$ is the maximal length of the paths from $X$ to $Y$. If $n_{\hspace{-0.5pt}_{XY}}=-1,$ then the claim follows easily from the above statement. If $n_{\hspace{-0.5pt}_{XY}}=0,$ then $\Hom_{\hspace{0.4pt}k\hspace{-0.5pt}\it\Delta}(X, Y)=k\hspace{0.4pt}\varepsilon_{\hspace{-0.5pt}_Y}$. On the other hand, $\Hom_{\hspace{1pt}\mathcal{C}}(X, Y)=k\id_{\hspace{-0.5pt}_Y}$ by the above sub-lemma, and the claim follows.

Suppose that $n_{\hspace{-0.5pt}_{XY}}>0.$ Let $\beta_i: Z_i\to Y$, $i=1, \ldots, s,$ be the arrows in $\Da$ ending in $Y$. Then $n_{\hspace{-0.5pt}_{XZ_i}}< n_{\hspace{-0.5pt}_{XY}}\vspace{0.5pt}$, and $(f_{\beta_1}, \cdots, f_{\beta_s}): Z_1\oplus \cdots \oplus Z_s\to Y$ is irreducible; see \cite[(3.4)]{Bau}. Since $Y$ admits a proper sink morphism, there exists a morphism $u: U\to Y$ such that $v=(f_{\beta_1}, \cdots, f_{\beta_s}, u): Z_1\oplus \cdots \oplus Z_s\oplus U\to Y$ is a proper sink morphism. Let $h: X\to Y$ be a morphism in $\C$. Being radical, $h$ factors through $v$. Since $\Da$ is a section, every indecomposable summand of $U$ lies in $\tau \Da;$ and since $\Hom_{\hspace{1pt}\mathcal{C}}(X, \Da^-)=0$, we have
$h=f_{\beta_1}h_1+ \cdots+f_{\beta_s}h_s$, with morphisms $h_i: X\to Z_i$ in $\C$. For each $1\le i\le s$, by the induction hypothesis, $h_i$ is a sum of composites of the chosen irreducible morphisms. Therefore, $h$ is a sum of composites of the chosen irreducible morphisms. Hence, $F_{_{XY}}$ is surjective.
Next, let $\rho: X\to Y$ be in $k\Da$ such that $F(\rho)=0$. Then $\rho=\beta_1\rho_1+\cdots+\beta_s\rho_s\vspace{1pt}$, where the $\rho_i: X\to Z_i$ are in $k\Da.$ Set
$w=(F(\rho_1), \cdots, F(\rho_s))^T: X\to Z_1\oplus \cdots\oplus Z_s.$ Then $(f_{\beta_1}, \cdots, f_{\beta_s})w=F(\rho)=0.$ If $\C$ has an almost slit sequence ending in $Y$
then, by Lemma \ref{auto-field}(2), $w$ factors through $\tau Y$; and since $\Hom_{\hspace{1pt}\mathcal{C}}(X, \Da^-)=0$, we have $w=0$. Otherwise, $v$ is a monomorphism, and hence, $w=0$. That is, in any case, $F(\rho_i)=0$, and by the inductive hypothesis, $\rho_i=0$, $i=1, \ldots, s$. As a consequence, $\rho=0$. Thus, $F_{_{XY}}$ is injective. This implies that $F$ is an equivalence. By Theorem \ref{mainprop}, $\Ga$ is standard. This establishes the sufficiency, and the necessity is evident.  The proof of the theorem is completed.

\medskip

Let $\Sa$ be a convex subquiver of $\GaC$. We shall say that $\Sa$ is {\it schurian} if, for any objects $X, Y$ in $\Sa$, the $k$-space $\Hom_{\mathcal{C}}(X, Y)$ is of dimension at most one; and it vanishes whenever $Y$ is a not successor of $X$ in $\Sa$. Moreover, we call $\Sa$ a {\it wing of rank $n$} if it is trivially valued of the following shape$\hspace{0.4pt}:$

\vspace{-8pt}
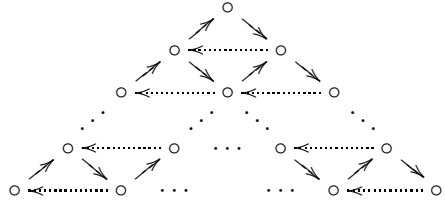
$$\xymatrixrowsep{16pt} \xymatrixcolsep{20pt}
\xymatrix@!=0.1pt{&&&& \circ \ar[dr]&&&&\\
&&&\circ \ar[ur]\ar[dr] \ar@{<.}[rr]&& \circ \ar[dr]&&&\\
&& \circ\ar[ur] \ar@{<.}[rr] && \circ\ar[ur] \ar@{<.}[rr] && \circ
&&\\
}$$ \vspace{-14pt}
$$\iddots \quad\;\; \quad \; \iddots \;\;\;  \ddots \;\quad \;\;\;\;\; \ddots$$
\vspace{-22pt}
$$\hspace{1pt}\xymatrixrowsep{16pt}
\xymatrixcolsep{20pt}\xymatrix@!=0.1pt{& \circ\ar[dr] \ar@{<.}[rr] && \circ   &\cdots& \circ\ar[dr]  \ar@{<.}[rr]  && \circ \ar[dr] & \\
\circ \ar[ur]\ar@{<.}[rr] && \circ \ar[ur]  &\cdots&&\cdots&
\circ\ar[ur] \ar@{<.}[rr]&& \circ,}\vspace{2pt}
$$
where the dotted arrows indicate the action of $\tau$, the objects are pairwise distinct and the number of $\tau$-orbits is $n\hspace{0.3pt};$ see \cite[(3.3)]{R3}. In this case, the object on the top is called the {\it wing vertex} and the objects at the bottom are said to be {\it quasi-simple}.

\medskip

\begin{Lemma}\label{wing}

Let $\mathcal{W}$ be a wing of $\GaC$. If the quasi-simple objects in $\mathcal{W}$ are pairwise orthogonal bricks, then $\mathcal{W}$ is schurian.

\end{Lemma}

\noindent {\it Proof.}  Assume that the quasi-simple objects in $\mathcal{W}$ are pairwise orthogonal bricks. Let $n$ be the rank of $\mathcal{W}$. If $n=1$, then the lemma holds trivially. Suppose that $n>1$ and the lemma holds for wings of rank $n-1$. Write the objects in $\mathcal{W}$ as $X_{ij}$ with $1\le j\le n$ and $j\le i \le n$ so that $X_{11}$ is the wing vertex and the $X_{nj}$ with $1\le j\le n$ are the quasi-simple objects. Observe that $X_{21}$ is the wing vertex of a schurian wing $\mathcal{W}_1$, while $X_{22}$ is the wing vertex of a schurian wing $\mathcal{W}_2$.
It is evident that we may choose irreducible morphisms $f_{ij}: X_{ij}\to X_{i+1, j}$ for $j\le i<n$ and $1\le j<n$; and irreducible morphisms $g_{pq}: X_{pq}\to X_{p-1,q-1}$ for $q\le p\le n$ and $2\le q\le n$ such that
\vspace{-8pt}
$$\qquad \mathcal{E}(X_{nj}): \; \xymatrix{X_{n,j+1}\ar[r]^-{g_{n,j+1}} & X_{n-1, j} \ar[r]^-{f_{n-1,j}}& X_{nj}}
$$ is an almost split sequence, for $j=1, \ldots, n-1;$ and

\vspace{-13pt}
$$\mathcal{E}(X_{ij}): \quad  \xymatrix{X_{i,j+1}\ar[rr]^-{(g_{i,j+1},\, f_{i,j+1})} && X_{i-1, j} \oplus X_{i+1,j+1}
\ar[rr]^-{\left(\hspace{-4pt}\begin{array}{c}f_{i-1,j}\\ g_{i+1,j+1}\end{array} \hspace{-5pt}\right)}
&& X_{i j}}$$
is an almost split sequence for $1\le j<n$ and $j<i<n$. Next, we shall divide the proof into several sub-lemmas.

\vspace{1pt}

(1) {\it $\Hom_{\hspace{0.4pt}\mathcal{C}}(X_{n1}, X_{ii})=0$ and $\Hom_{\,\mathcal{C}}(X_{i1}, X_{nn})=0$, for $1\le i\le n$}. Suppose that $\C$ has a non-zero morphism $f: X_{n1}\to X_{rr}$ for some $1\le r\le n$. Assume that $r$ is maximal.  Since $X_{n1}, X_{nn}$ are orthogonal, we have $r<n$. Since $\mathcal{W}_1$ is schurian, $f_{rr}f=0$. Applying Lemma \ref{auto-field}(2) to the almost split sequence $\mathcal{E}(X_{r+1,r})$, we see that $f$ factors through $g_{r+1,r+1}: X_{r+1, r+1}\to X_{rr}$, which contradicts the maximality of $r$. The first part of the statement is established. In a dual manner, we may prove the second part.

\vspace{1pt}

(2) {\it $\Hom_{\hspace{0.4pt}\mathcal{C}}(X_{i1}, \mathcal{W}_2)=0$ and $\Hom_{\hspace{0.4pt}\mathcal{C}}(\mathcal{W}_1, X_{ii})=0,$ for $1\le i\le n$.} Suppose that $f: X_{s1}\to X$ is a non-zero morphism with $1\le s\le n$ and $X\in \mathcal{W}_2$, which is necessarily  radical. If $X\ne X_{jj}$ for any $2\le j\le n$, then $X$ admits a sink morphism whose domain is a direct sum of one or two objects in $\mathcal{W}_2$. Factorizing $f$ through this sink morphism, we obtain a non-zero morphism $g: X_{s1}\to X_{ii}$ with $2\le i\le n$. Assume that $s$ is maximal for this property. By Statement (1), $s<n$. Since $\mathcal{W}_2$ is schurian, $gg_{s+1,2}=0$. Applying Lemma \ref{auto-field}(3) to $\mathcal{E}(X_{s+1,1})$, we see that $g$ factors through $f_{s1}: X_{s1}\to X_{s+1,1}$, which contradicts the maximality of $s$. The first part of the statement is established. In a dual fashion, we may establish the second part.

\vspace{1pt}

(3) {\it $\Hom_{\,\mathcal{C}}(X_{nn}, \mathcal{W}_1)=0$ and $\Hom_{\,\mathcal{C}}(\mathcal{W}_2, X_{n1})=0.$} Suppose that $\mathcal{C}$ has a non-zero morphism $f: X_{nn}\to X_{pq}$ with $2\le p\le n$ and $1\le q <p$. We may assume that $p$ is maximal for this property. Since the quasi-simple objects are orthogonal, $p<n$. By the maximality of $p$, we have $f_{pq}f=0$. Applying Lemma \ref{auto-field}(2) to $\mathcal{E}(X_{p+1,q})$, we see that $f$ factors through $g_{p+1, q+1}$, contrary to the maximality of $p$. The first part of the statement is established, and the second part follows dually.

\vspace{1pt}

(4) {\it If $f: X_{ii}\to X_{11}$ with $1\le i<n$ is such that $f g_{i+1, i+1}\cdots g_{nn}=0$, then $f=0$. Dually, if $g: X_{11}\to X_{i1}$ is a morphism with $1\le i<n$ such that $f_{n-1,1}\cdots f_{i1} g=0$, then $g=0$.} Suppose that $f g_{i+1, i+1}\cdots g_{nn}=0$ but $f\ne 0$. Let $r$ with $i+1\le r\le n$ be minimal such that $f g_{i+1, i+1}\cdots g_{rr}=0$. Write $f g_{i+1, i+1}\cdots g_{rr}=g g_{rr}$, where $g: X_{r-1,r-1}\to X_{11}$ is a non-zero morphism. Applying Lemma \ref{auto-field}(3) to $\mathcal{E}(X_{r, r-1})$, we see that $g$ factors through $f_{r-1, r-1}$, which contradicts Statement (2). This establishes the first par of the statement.

\vspace{1pt}

%

(5)  {\it $\Hom_{\hspace{0.4pt}\mathcal{C}}(X_{ii}, X_{11})$ and $\Hom_{\hspace{0.4pt}\mathcal{C}}(X_{11}, X_{i1})$ are one-dimensional, for $1\le i\le n$.} It suffices to prove the first part of the statement, since the second part follows dually. Let $f: X_{nn}\to X_{11}$ be a morphism. By Statement (3), $f_{11} f=0$. Applying Lemma \ref{auto-field}(2) to $\mathcal{E}(X_{21}),$ we obtain some $f_1: X_{nn}\to X_{22}$ such that $f=g_{22}f_1$. Since $f_{22}f_1=0$ by Statement (3), we may repeat this process to obtain a morphism $f_{n-1}: X_{nn} \to X_{nn}$ such that $f=g_{22}\cdots g_{nn} f_{n-1}$. Since $\mathcal{W}_1$ is schurian, $f_{n-1}=\lambda \id$ for some $\lambda\in k$, and hence, $f=\lambda  g_{22}\cdots g_{nn}$. Since $g_{22}\cdots g_{nn}\ne 0$, we see that $\{g_{22}\cdots g_{nn}\}$ is a $k$-basis for $\Hom_{\hspace{0.4pt}\mathcal{C}}(X_{nn}, X_{11})$. Write $g_{11}=\id_{X_{11}}$. If $g: X_{ii}\to X_{11}$ is a morphism with $1\le i<n$, then $gg_{i+1,i+1} \cdots g_{nn}=\mu g_{22}\cdots g_{nn}=\mu g_{11}\cdots g_{nn}$, for some $\mu\in k$. This yields that
$(g-\mu g_{11}\cdots g_{ii}) g_{i+1, i+1}\cdots g_{nn}=0.$ By the first part of Statement (4), $g=\mu g_{11}\cdots g_{ii}.$ Being non-zero,  $g_{11}\cdots g_{ii}$ forms a $k$-basis for $\Hom_{\hspace{0.4pt}\mathcal{C}}(X_{ii}, X_{11})$.

\vspace{1pt}

Now, suppose that $\Hom_{\hspace{0.4pt}\mathcal{C}}(X, Y)\ne 0$ for some $X, Y\in\mathcal{W}.$ We claim that $Y$ is a successor of $X$ and $\Hom_{\hspace{0.4pt}\mathcal{C}}(X, Y)$ is one-dimensional. If $X\in \mathcal{W}_1$, then $Y\in \mathcal{W}_1$ by Statement (2). Since $\mathcal{W}_1$ is schurian, our claim follows. Otherwise, $X=X_{ss}$ for some $1\le s\le n$. If $s=n$ then, by Statement (3), $Y=X_{ii}$ for some $1\le i\le n$. Combining Statement (5) and the fact that $\mathcal{W}_2$ is schurian, we see that $\Hom_{\hspace{0.4pt}\mathcal{C}}(X, Y)$ is one-dimensional. If $s=1$, then $Y=X_{j1}$ for some $1\le j\le n$, and hence, $\Hom_{\hspace{0.4pt}\mathcal{C}}(X, Y)$ is one-dimensional by Statement (5).
Finally, suppose that $1<s<n$. If $Y\in \mathcal{W}_2$, since $\mathcal{W}_2$ is schurian, our claim follows. Otherwise, by Statement (3), $Y=X_{t1}$ for some $1\le t< n$. If $t=1$, then $\Hom_{\hspace{0.4pt}\mathcal{C}}(X, Y)$ is one-dimensional by Statement (5). It remains to consider the case where $1<t<n$. Let $f: X_{ss}\to X_{t1}$ be a non-zero morphism with $1<s, t<n$. Factorizing $f$ along the $\mathcal{E}(X_{j1})$ with $2\le j\le t$, we get $g: X_{ss}\to X_{t+1, 2}$ and $h: X_{ss}\to X_{11}$ such that $f=g_{t+1, 2}g+f_{t-1,1}\cdots f_{11} h.$  By Statement (5), $h=\lambda g_{22}\cdots g_{ss}$ with $\lambda\in k$. This yields $f=g_{t+1, 2} u$, where $u: X_{ss}\to X_{t+1, 2}$ is a non-zero morphism. Since $\mathcal{W}_2$ is schurian, $X_{t+1, 2}$ is a successor of $X_{ss}$ and $\Hom_{\hspace{0.4pt}\mathcal{C}}(X_{ss}, X_{t+1,2})$ has a $k$-basis $\{v\}$. Therefore, $f=\mu g_{t+1, 2} v$ with $\mu\in k$. This shows that $\{g_{t+1, 2} v\}$ is a $k$-basis for $\Hom_{\hspace{0.4pt}\mathcal{C}}(X_{ss}, X_{t1})$. This establishes our claim. The proof of the lemma is completed.

\medskip

Let $\mathbb{A}_{\tiny\infty}^+$ and $\mathbb{A}_{\tiny\infty}^-$ denote the linearly oriented quivers of type $\A_\infty$ having a unique source and having a unique sink, respectively. If $\Ga$ is a connected component of $\GaC$ of shape $\mathbb{Z}\mathbb{A}_{\tiny\infty}$, $\mathbb{N}\mathbb{A}_{\tiny\infty}^+$ or $\mathbb{N}^-\mathbb{A}_{\infty}^-,$ then the objects in $\Ga$ having at most one immediate predecessor and at most one immediate successor are called {\it quasi-simple.}

\medskip

\begin{Theo}\label{MainSemiStable}

Let $\C$ be a Hom-finite Krull-Schmidt additive $k$-category with $\Ga$ a connected component of $\GaC$. If $\Ga$ is a wing or of shape  $\mathbb{Z}\mathbb{A}_{\tiny\infty}$, $\mathbb{N}\mathbb{A}_{\tiny\infty}^+$ or $\mathbb{N}^-\mathbb{A}_{\infty}^-,$ then it is standard if and only if its quasi-simple objects are pairwise orthogonal bricks.

\end{Theo}

\noindent{\it Proof.} 
We shall need only to prove the sufficiency. Let $\Ga$ be a wing or of shape $\mathbb{Z}\mathbb{A}_{\tiny\infty}$, $\mathbb{N}\mathbb{A}_{\tiny\infty}^+$ or $\mathbb{N}^-\mathbb{A}_{\infty}^-$ with the quasi-simple objects being pairwise orthogonal bricks. Then any two objects in $\Ga$ lie in a wing whose quasi-simples are pairwise orthogonal bricks. By Lemma \ref{wing}, $\Ga$ is schurian. Choose a section $\Da$ of $\Ga$ so that $\Da$ is the right-most section if $\Ga$ is a wing or of shape $\mathbb{N}^-\mathbb{A}_{\infty}^-;$  and $\Da$ is the left-most section if $\Ga$ is of shape $\N\A^+_\infty;$ and $\Da$ is any section with an alternating orientation if $\Ga$ is of shape $\Z\A_\infty$. Then $\Da^-$ has no right infinite path and $\Da^+$ has no left infinite path such that $\Hom_{\hspace{0.5pt}\mathcal{C}}(\Da^+, \Da\cup\Da^-)=0$ and $\Hom_{\hspace{0.5pt}\mathcal{C}}(\Da, \Da^-) = 0.$

For each arrow $\alpha: X\to Y$ in $\Da$, we choose an irreducible morphism $f_\alpha: X\to Y$ in $\C$. Since every path in $\Da$ is sectional, the composite of any chain of the chosen irreducible morphisms is non-zero; see \cite[(2.7)]{L3}. Therefore, for any $M, N\in \Da$, $\Hom_{\,\mathcal{C}}(M, N)$ is one-dimensional if and only if $N$ is a successor of $M$ in $\Da$, and in this case, the composite of the chain of the chosen irreducible morphisms corresponding to the path from $M$ to $N$ forms a $k$-basis for $\Hom_{\,\mathcal{C}}(X, Y).$ It is now easy to see that $k\Da \cong \C(\Da).$ By Theorem \ref{mainprop}, $\Ga$ is standard. This establishes the sufficiency, and the necessity is trivial. The proof of the theorem is completed.

\section{Specialization to representation categories of quivers}

Throughout this section, we fix a connected quiver $Q=(Q_0, Q_1)$, where $Q_0$ is the set of vertices and $Q_1$ is the set of arrows, which is assumed to be strongly locally finite, that is, $Q$ is locally finite such that
the number of paths between any two given vertices is finite. A {\it $k$-representation} $M$ of $Q$ consists of a family of $k$-spaces $M(x)$ with $x\in Q_0,$ and a family of $k$-maps $M(\alpha): M(x)\to M(y)$ with $\alpha: x\to y\in Q_1$. For such a representation $M$, one defines its {\it support} ${\rm supp}\,M$ to be the full subquiver of $Q$ generated by the vertices $x$ for which $M(x)\ne 0$, and one calls $M$
\emph{locally finite dimensional} if ${\rm dim}_k M(x)\vspace{0pt}$ is finite for all $x\in Q_0,$ and
{\it finite dimensional} if $\Sigma_{x \in Q_0}{\rm dim}_k M(x)$ is finite. The locally finite dimensional $k$-representations of $Q$ form a hereditary abelian $k$-category $\rep(Q)$. The subcategory of $\rep(Q)$ of finite dimensional representations is written as $\rep^b(Q)$. For each $x\in Q_0$, one constructs an indecomposable projective representation $P_x$ and an indecomposable injective representation $I_x$; see \cite[Section 1]{BLP}. Since $Q$ is strongly locally finite, $P_x$ and  $I_x$ lie in $\rep(Q)$. One says that $M\in \rep(Q)$ is {\it finitely presented} if $M$ has a minimal projective presentation $\xymatrixcolsep{14pt}\xymatrix{P_1 \ar[r]&
P_0 \ar[r] & M \ar[r]& 0,\hspace{-0pt}}$ where $P_1, P_0$ are finite direct sums of some $P_x$ with $x\in Q_0$; and {\it finitely co-presented} if $M$ has a minimal injective co-presentation $\xymatrixcolsep{14pt}\xymatrix{0\ar[r]& M \ar[r] &
I_0 \ar[r] & I_1,}\hspace{-2pt}$ where $I_0, I_1$ are finite direct sums of some $I_x$ with $x\in Q_0$.
Let $\rep^+(Q)$ and $\rep^-(Q)$ be the full subcategories of $\rep(Q)$ of finitely presented representations and of finitely co-presented representations, respectively. Then $\rep^b(Q)$ is the intersection of $\rep^+(Q)$ and $\rep^-(Q)$. In particular, $I_x\in \rep^+(Q)$ if and only if $I_x\in \rep^b(Q)$. One denotes by $Q^+$ the full subquiver of $Q$ generated by the vertices $x$ for which $I_x\in \rep^b(Q).$

\medskip

It is known that $\rep^+(Q)$ and $\rep^-(Q)$ are hereditary, abelian and Hom-finite; see \cite[(1.15)]{BLP}. In particular, they are Krull-Schmidt.
The shapes of the their Auslander-Reiten components have been well described.
Indeed, the  Auslander-Reiten quiver $\Ga_{\rep^+(Q)}\vspace{1pt}$ of $\rep^+(Q)$ has a unique {\it preprojective component}, which has a left-most section generated by the $P_x$ with $x \in Q_0$; see \cite[(4.6)]{BLP} and \cite[(2.4)]{R3}. The connected components of $\Ga_{\rep^+(Q)}\vspace{1pt}$ containing some of the $I_x$ with $x\in Q^+$ are called {\it preinjective}, which correspond bijectively to the connected components of the quiver $Q^+$. Note that every preinjective component has a right-most section generated by its injective representations $I_x$; see \cite[(2.4)]{R3} and \cite[(4.7)]{BLP}. The other connected components of $\Ga_{\rep^+(Q)}$ are called \emph{regular}, which are wings, stable tubes or of shapes $\Z\A_\infty, \N \A^+_\infty$ and $\N^-\A^-_\infty;$ see \cite[(4.14)]{BLP}, \cite{Rin} and \cite{R3}.

\medskip

The following easy fact is well known in the finite case.

\medskip

\begin{Lemma} \label{classicalAR}

Let $X$ and $Y$ be representations lying in $\Ga_{\rep^+(Q)}\vspace{1.4pt}$. If $\tau X$ and $\tau Y$ are defined in $\Ga_{\rep^+(Q)}$, then $\Hom_{\rep^+(Q)}(X, Y)\cong \Hom_{\rep^+(Q)}(\tau X, \tau Y).$

\end{Lemma}

\noindent{\it Proof.} Assume that  $\tau X$ and $\tau Y$ are defined in $\Ga_{\rep^+(Q)}\vspace{1pt}$. In view of the proof stated in \cite[(2.8)]{BLP}, we have $\Hom(\tau X, \tau Y)\cong D\Ext^1(Y, \tau X)$.  Dually, since $\tau X$ is not injective  and finite-dimensional; see \cite[(3.6)]{BLP}, \vspace{0pt} $\Hom(X, Y)\cong D\Ext^1(Y, \tau X)$. The proof of the lemma is completed.

\medskip

Recall that $Q$ is of \emph{infinite Dynkin type} if its
underlying graph is $\A_\infty$, $\A_\infty^\infty$ or $\mathbb{D}_\infty.$
In this case, a reduced walk is a {\it string} if it contains at most finitely many, but at least one, sinks or sources. To each  string $w$, one associates a {\it string representation} $M_w$ defined as follows: for $x\in Q_0$, one sets $M_w(x)=k$ if $x$ appears in $w$, and otherwise, $M_w(x)=0$; and for $\alpha\in Q_1$, one sets $M_w(\alpha)=\id$ if $\alpha$ arrears in $w$, and otherwise, $M_w(\alpha)=0;$ see \cite[Section 5]{BLP}. It is easy to see that every string representation has a trivial endomorphism algebra.

\medskip

\begin{Theo} \label{standardpreproj}

Let $Q$ be a connected quiver which is strongly locally finite.

\vspace{-1.5pt}

\begin{enumerate}[$(1)$]

\item The preprojective component and the preinjective components of $\Ga_{\rep^+(Q)}$ are standard.

\item If $Q$ is of finite or infinite Dynkin type, then every connected component of $\Ga_{\rep^+(Q)}$ is standard.

\vspace{1pt}

\item If $Q$ is infinite but not of infinite Dynkin type, then $\Ga_{\rep^+(Q)}$ has infinitely many non-standard regular components.

\end{enumerate}

\end{Theo}

\noindent{\it Proof.} (1) 
The preprojective component $\mathcal{P}_Q$ of $\Ga_{\rep^+(Q)}\vspace{1pt}$ has a left-most section
$\Da$ which is generated by the $P_x$ with $x \in Q_0$ and isomorphic to $Q^{\, \rm op}$; see \cite[(4.6)]{BLP} and \cite[(2.4)]{R3}. Hence, $\Da^- = \emptyset$. Moreover, $\Da^+$ has no left infinite path; see \cite[(4.8)]{BLP} and (\ref{no-inf-path}). If $f: X\to Y$ is a non-zero morphism with $X\in \mathcal{P}_Q$ and $Y\in \Da$, then $X$ is a predecessor of $Y$ in $\mathcal{P}_Q$; see \cite[(4.9)]{BLP}, and hence, $X\in \Da$. Therefore, $\Hom_{\rep^+(Q)}(\Da^+, \Da)=0\vspace{1pt}$. Let $\mathcal{P}(Q)$ be the full subcategory of $\rep^+(Q)$ generated by the $P_x$ with $x\in Q_0$. For each arrow $\alpha: y\to x$ in $Q$, denote by $P_\alpha: P_x\to P_y$ the morphism given by the right multiplication by $\alpha$. It is easy to see that
$$F: kQ^{\rm \,op} \to \mathcal{P}(Q): x\mapsto P_x; \hspace{3pt} \alpha^{\rm o}\mapsto P_\alpha \vspace{-1pt}$$
is a faithful $k$-functor, which is also full by Proposition 1.3 stated in \cite{BLP}. Thus, $\Da$ is standard.
By Theorem \ref{mainprop}, $\mathcal{P}_Q$ is standard. Dually, the preinjective component $\mathcal{I}\hspace{0.4pt}_Q$ of $\Ga_{\rep^-(Q)}\vspace{0.5pt}$ is standard. By the dual of Lemma 4.5(1) stated in \cite{BLP}, $\mathcal{I}\hspace{0.4pt}_Q$ has a left-most section $\mathit\Theta$ generated by its infinite-dimensional representations. Now the preinjective components of $\Ga_{\rep^+(Q)}\vspace{0.5pt}$ are the connected components of the complement of  $\it\Theta$ in $\mathcal{I}\hspace{0.4pt}_Q$. Hence, each preinjective component of $\Ga_{\rep^+(Q)}$ is a convex translation subquiver of $\mathcal{I}_Q$, and in particular, it is standard.

(2) Suppose that $Q$ is of infinite Dynkin type. Let $\Ga$ be a regular component of $\Ga_{\rep^+(Q)}\vspace{0.5pt}.$ Then $\Ga$ is a wing or of shape $\Z A_\infty$, $\N^-\A^-_\infty$ or $\N\A_\infty^+$; see \cite[(4.14)]{BLP}. Moreover, $Q$ is of type $\A_\infty^\infty\vspace{0.5pt}$ or $\mathbb{D}_\infty$; see \cite[(5.16)]{BLP}.
Assume first that $Q$ is of type $\mathbb{A}_{\infty}^{\infty}$. By Proposition 5.9 stated in \cite{BLP}, the representations in $\Ga_{\rep^+(Q)}\vspace{1pt}$ are all string representations, and hence, they are all bricks. Moreover, the quasi-simple representations in $\Ga$ have pairwise disjoint supports; see \cite[(5.15)]{BLP}. In particular, they are pairwise orthogonal. By Theorem \ref{MainSemiStable}, $\Ga$ is standard.

Assume next that $Q$ is of type $\mathbb{D}_{\infty}$. Then $\Ga$
is of shape $\mathbb{Z}\mathbb{A}_{\infty}$, $\mathbb{N}\mathbb{A}^+_\infty$ or $\mathbb{N}^-\mathbb{A}^-_{\infty};$ see \cite[(5.22)]{BLP}. In particular, $\tau$ or $\tau^-$ is defined everywhere in $\Ga$. We shall consider only the first case, since the second case can treated in a dual manner. Let $a\in Q_0$ be one of the two vertices of weight one, which lies in the support of at most two quasi-simple representations; see \cite[(5.20)]{BLP}. Thus, there exists a quasi-simple representation $S\in \Ga$ such that $(\tau^nS)(a)=0$, for all $n \ge 0$.

Let $M, N$ be quasi-simple representations in $\Ga$. There exists $m\ge 0$ such that $\tau^mM=\tau^rS$ and $\tau^mN=\tau^sS$ with $r, s\ge 0$. We may assume that $r\ge s$. By Lemma \ref{classicalAR},
$\Hom(M, N) \cong \Hom(\tau^mM, \tau^mN) = \Hom(\tau^rS,\tau^sS).$ Since $(\tau^rS)(a)=0$, we see that $\tau^rS$ is a string representation; see \cite[(5.19)]{BLP}. Thus, $\tau^rS$ is a brick. Taking $N=M$, we see that $M$ is a brick. Suppose that $M\ne N$. Then $r>s$. Setting $t=r-s$, there exists a sectional path $\xymatrixcolsep{16pt}\xymatrix{\hspace{-2pt}S_t\ar[r] & S_{t-1}\ar[r] & \cdots \ar[r] & S_1 \ar[r] & \tau^sS}$ in $\Ga.$ For $x\in Q_0$, we have ${\rm dim} \hspace{0.3pt}S_t(x) = {\textstyle\sum}_{i=s}^{\hspace{0.3pt}r}\,{\rm dim}\hspace{0.3pt}\tau^iS(x)$. Since $(\tau^iS)(a)=0$ for $i\ge 0,$ ${\rm dim}\hspace{0.3pt}S_t(a)=0$. Hence, $S_t$ is a string representation; see \cite[(2.19)]{BLP}. If the supports of $\tau^rS$ and $\tau^sS$ have a common vertex $b$, then ${\rm dim} \hspace{0.3pt}S_t(b)\ge {\rm dim}\hspace{0.3pt}\tau^rS(b) + {\rm dim}\hspace{0.3pt}\tau^sS(b)\ge 2$, contrary to $S_t$ being a string representation. Thus, $\tau^rS$ and $\tau^sS$ have disjoint supports. In particular, they are orthogonal, and so are $M$ and $N$. By Theorem \ref{MainSemiStable}, $\Ga$ is standard. In view of Statement (1), we have established Statement (2).

(3) Suppose that $Q$ is infinite but not of infinite Dynkin type. Then $Q$ has a finite subquiver $\Sa$ of Euclidean type. Then we can find a homogeneous tube $\mathcal{T}$ in $\Ga_{\rep^b(\it\Sigma)}$; see, for example, \cite[(6.3)]{BLP}. Let $M_i$ with $i\ge 1$ be the representations in $\mathcal{T}$ which are not quasi-simple. Regarded as representations of $Q$, the $M_i$ are distributed in infinitely many regular components of $\Ga_{\rep^+(Q)}$; see \cite[(6.1), (6.2)]{BLP}. These regular components are not standard, since the $M_i$ have non-trivial endomorphism algebras. The proof of the theorem is completed.

\medskip

\noindent{\sc Remark.} (1) In view of Theorem 5.17 stated in \cite{BLP}, we see that wings and the translation quivers $\Z\A_\infty$, $\N \A^+_\infty$ and $\N^-\A^-_\infty$ all occur as standard Auslander-Reiten components of Krull-Schmidt categories.

\vspace{2pt}

\noindent (2) Let $Q$ be finite of Euclidean type. If $k$ is not algebraically closed, then some indecomposable $k$-representations of $Q$ have a non-trivial automorphism field; see the proof in \cite[(6.3)]{BLP}. As a consequence, every connected component of $\Ga_{\rep^b(Q)}$ is standard if and only if $k$ is algebraically closed.

\medskip

We conclude this section with an application to the bounded derived category $D^b(\rep^+(Q))$ of $\rep^+(Q)$. Since  $\rep^+(Q)$ is hereditary, the vertices of $\Ga_{D^b(\rep^+(Q))}\vspace{1pt}$ can be chosen to be the shifts of those in $\Ga_{\rep^+(Q)}\vspace{1pt}$. If $Q$ is not of finite Dynkin type, then the connected components of $\Ga_{D^b(\rep^+(Q))}\vspace{0.5pt}$ are the shifts of the regular components of $\Ga_{\rep^+(Q)}$ and the shifts of the {\it connecting component}, which is obtained by gluing the preprojective component together with the shift by $-1$ of the preinjective components of $\Ga_{\rep^+(Q)};$ see \cite[(4.4)]{Hap} and \cite[(7.10)]{BLP}.
In case $Q$ is of finite Dynkin type, $\Ga_{D^b(\rep^+(Q))}\vspace{1pt}$ is connected of shape $\Z Q^{\hspace{0.4pt}\rm op}$, which is obtained by gluing, for each integer $i$, the shift by $i$ of $\Ga_{\rep^b(Q)}$ together with its shift by $i+1;$
see \cite[(4.5)]{Hap}. In this case, we also call $\Ga_{D^b(\rep^+(Q))}\vspace{1pt}$ the {\it connecting component}.

\smallskip

\begin{Theo} \label{connecting} If $Q$ is a connected strongly locally finite quiver, then the connecting component of $\Ga_{D^b(\rep^+(Q))}\vspace{1pt}$ is standard$\hspace{0.3pt};$ and every connected component is standard in case $Q$ is of finite or infinite Dynkin type.
\end{Theo}

\noindent{\it Proof.} Assume that $Q$ is a connected strongly locally finite quiver and $\mathcal{C}_Q$ is the connecting component of $\Ga_{D^b(\rep^+(Q))}\vspace{1pt}$. Let $\Da$ be the full subquiver of $\C_Q$ generated by the representations $P_x\in \Ga_{\rep^+(Q)}$ with $x\in Q_0$, which is isomorphic to $Q^{\rm op}$. It follows from Lemma 7.8 stated in \cite{BLP} that $\Da$ is a section of $\C_Q$. Since $\rep^+(Q)$ fully embeds in $D^b(\rep^+(Q))$, by Theorem \ref{standardpreproj}, $\Da$ is standard. Let $M, N\in \rep^+(Q)$. Since $\rep^+(Q)$ is hereditary, $\Hom_{D^b(\rep^+(Q))}(M[m], N[n])=0$ for $m>n$; see \cite[(3.1)]{Lez}. Combining this fact with the  standardness of the preprojective component of $\Ga_{\rep^+(Q)}$, we deduce easily that $\Hom_{D^b(\rep^+(Q))}(\Da^+, \Da)=0$ and $\Hom_{D^b(\rep^+(Q))}(\Da\cup\Da^+, \Da^-)=0.\vspace{1pt}$

If $Q$ is not of finite Dynkin type, then $\Da^+$ coincides with the full subquiver of the preprojective component of $\Ga_{\rep^+(Q)}$ generated by the non-projective representations, while $\Da^-$ coincides with the shift by $-1$ of the preinjective components of $\Ga_{\rep^+(Q)}$. Thus, $\Da^+$ contains no left infinite path and $\Da^-$ contains no right infinite path by Lemma 4.8 stated in \cite{BLP}. This is also the case if $Q$ is of finite Dynkin type; see  (\ref{no-inf-path}). Thus $\C_Q$ is standard by Theorem \ref{mainprop}. This establishes the first part of the theorem. Combining this with Theorem \ref{standardpreproj}(2), we obtain the second part. The proof of the theorem is completed.

\medskip

\noindent{\sc Remark.} Let $Q$ have no infinite path. If $Q$ is not of finite Dynkin type, then $\Ga_{\rep^+(Q)}$ has a unique preinjective component of shape $\N Q^{\hspace{0.4pt}\rm op}$ and its proprojective component is of shape $\N^-\hspace{-0.5pt}Q^{\hspace{0.4pt}\rm op}$; see \cite[(4.7)]{BLP}. Thus, in any case, the connecting component of $\Ga_{D^b(\rep^+(Q))}$ is standard of shape $\Z Q^{\hspace{0.4pt}\rm op}$.

\section{Specialization to module categories of algebras}

Throughout this section, assume that $k$ is algebraically closed. Let $A$ stand for a finite-dimensional $k$-algebra and ${\rm mod}\hspace{0.4pt}A$ for the category of finite-dimensional left $A$-modules. In this classical situation, we have the following easy criteria for an Auslander-Reiten component with sections to be standard.

\medskip

\begin{Theo} \label{maincor}

Let $A$ be a finite-dimensional algebra over an algebraically closed field, and let $\Ga$ be a connected component of
$\Ga_{{\rm mod}\hspace{0.3pt}A}$. If $\Da$ is a section of $\Ga$, then $\Ga$ is standard if and only if
$\Hom_A(\Da, \tau\hspace{-1pt}\Da)=0$ if and only if $\Hom_A(\tau^-\hspace{-1pt}\Da, \Da)=0.$

\end{Theo}

\noindent{\it Proof.} Let $\Da$ be a section of $\Ga$. Note that every module in $\Da$ admits a proper sink map and a proper source map. Moreover, since the base field is algebraically closed, every module in $\Da$ has a trivial automorphism field.

Suppose that $\Hom_A(\Da, \tau \Da)=0.$ Then $\Da$ is finite; see \cite[(2.1)]{Skow}. By Lemma \ref{no-inf-path}, $\Da^+$ has no left infinite path and $\Da^-$ has no right infinite path. Assume that $\Hom_A(X, Y)\ne 0$ for some $X\in \Da^+$ and $Y\in \Da^-$. Since every module in $\Da^+$ admits a sink epimorphism, we obtain an arrow $X_1\to X$ in $\Ga$ such that $\Hom_A(X_1, Y)\ne 0$. Observe that $X_1\in \Da\cup \Da^+$. If $X_1\in \Da^+$, then $\Ga$ has an arrow $X_2\to X_1$ such that $\Hom_A(X_2, Y)\ne 0$. Since $\Da^+$ has no left infinite path, there exists a module $M$ in $\Da$ such that $\Hom_A(M, Y)\ne 0.$ Similarly, since $\Da^-$ has no right infinite path and every module in $\Da^-$ has a source monomorphism, there exists a module $N$ in $\tau \Da$ such that $\Hom_A(M, N)\ne 0,$ a contradiction. This shows that $\,\Hom_A(\Da^+, \Da^-)=0.$ By Theorem \ref{Theorem2}, $\Ga$ is standard. If $\Hom_A(\tau^-\Da, \Da)=0,$ one shows in a dual manner that $\Ga$ is standard. Conversely, it is evident that $\Hom_A(\Da, \tau \Da)=0$ and $\Hom_A(\tau^-\Da, \Da)=0$ if $\Ga$ is standard. The proof of the theorem is completed.

\medskip

Let $\Ga$ be a connected component of $\Ga_{{\rm mod}\hspace{0.3pt}A}$. Recall that $\Ga$ is {\it generalized standard} if ${\rm rad}^\infty({\rm mod}\hspace{0.4pt}A)$ vanishes in $\Ga$; see \cite{Skow}. It is known that $\Ga$ is generalized standard if it is standard; see \cite{L6}, and the converse holds true in case $\Ga$ has no projective module or no injective module; see \cite{Skow3}. Observing that the conditions on $\Da$ stated in Theorem \ref{maincor} are trivially verified in case $\Ga$ is generalized standard, we obtain the following consequence.

\medskip

\begin{Cor} \label{cor1}
Let $\Ga$ be a connected component of $\Ga_{{\rm mod}\hspace{0.4pt}A}$. If $\Ga$ has a section, then it is standard if and only if it is generalized standard.
\end{Cor}

\medskip

The algebra $A$ is called {\it tilted} if $A=\End_H(T)$, where $H$ is a finite-dimensional hereditary algebra and $T$ is a tilting $H$-module. In this case, ${\rm mod}\hspace{0.3pt}A$ contains slices; see \cite{HaR}, and a connected component of $\Ga_{{\rm mod}\hspace{0.3pt}A}$ containing the indecomposable mo\-dules of a slice is called a {\it connecting component}\hspace{0.3pt}. It is shown that a connecting component of a tilted algebra is standard; see \cite[(5.7)]{ABM}.

\medskip

\begin{Cor} \label{cor2}

If $\Ga$ is a connected component of $\Ga_{{\rm mod}\hspace{0.4pt}A}$, then $\Ga$ is standard with sections if and only if it is a connecting component of a tilted factor algebra of $A$.

\end{Cor}

\noindent{\it Proof.} Let $\Ga$ be a connected component of $\Ga_{{\rm mod}\hspace{0.4pt}A}$. Suppose first that $\Ga$ is standard with a section $\Da$. In particular, we have ${\rm Hom}_A(\Da, \tau \Da)=0.$ If $I$ is the intersection of the annihilators of the modules in $\Ga$, then $B=A/I$ is a tilted algebra with $\Ga$ a connecting component of $\Ga_{{\rm mod}\hspace{0.4pt}B};$ see \cite[(2.2)]{L5}, and also \cite{Skow2}.

Suppose next that there exists a tilted algebra $B=A/I$ with $\Ga$ being a connecting component of $\Ga_{{\rm mod}\hspace{0.4pt}B}$. Then $\Ga$ has a section $\Da$ generated by the non-isomorphic indecomposable modules of a slice of ${\rm mod}\hspace{0.4pt}B.$ By the defining property of a slice, ${\rm Hom}_B(\Da, \tau \Da)=0.$ Thus, by Theorem \ref{maincor}, $\Ga$ is a standard component of $\Ga_{{\rm mod}\hspace{0.4pt}B}$. Since ${{\rm mod}\hspace{0.4pt}B}$ fully embeds in ${{\rm mod}\hspace{0.4pt}A}$, we see that $\Ga$ is a standard component of $\Ga_{{\rm mod}\hspace{0.4pt}A}$. The proof of the corollary is completed.

\bigskip\bigskip

\small
$\begin{array}{lllllll}
\mbox{Shiping Liu}                       && \mbox{Charles Paquette} \\
\mbox{D\'epartement de math\'ematiques}  && \mbox{Department of Mathematics and Statistics} \\
\mbox{Universit\'e de Sherbrooke}        && \mbox{University of New Brunswick} \\
\mbox{Sherbrooke (QC), Canada J1K 2R1}   && \mbox{Fredericton (NB), Canada E3B 9P8}\\
\mbox{Email: shiping.liu@usherbrooke.ca} && \mbox{Email: charles.paquette@usherbrooke.ca}
\end{array}$

\end{document}